\documentclass[10pt,oneside]{amsart}

\usepackage{amssymb, amscd, amsmath, amsthm, latexsym}
\usepackage{pb-diagram, fancyhdr, graphicx, psfrag, color}
  
\newcommand{\compactlist}{\begin{list}{$\bullet$}{\setlength{\leftmargin}{1em}}}
\def\zz{{\bf Z}}

\def\ff{{\bf F}}
 
\def\qq{{\bf Q}}

\def\rr{{\bf R}}

\def\calf{\mathcal{F}}

\def\calc{\mathcal{C}}

\def\calg{\mathcal{G}}

\newcommand{\spinc}{\ifmmode{{\mathfrak s}}\else{${\mathfrak s}$\ }\fi}
\newcommand{\spinct}{\ifmmode{{\mathfrak t}}\else{${\mathfrak t}$\ }\fi}

\newtheorem{theorem}{Theorem}
\newtheorem{lemma}[theorem]{Lemma}

\theoremstyle{definition}

\numberwithin{equation}{section}

\begin{document}
\title[The concordance classification of knots]{The concordance classification of low crossing number knots}
\author{Julia Collins}\author{Paul Kirk}\author{Charles Livingston}

\thanks{This work was supported in part by the National Science Foundation under Grant  1007196. \\ \today}

\address{Julia Collins: Department of Mathematics, University of Edinburgh, Scotland}
\email{ Julia Collins <Julia.Collins@ed.ac.uk>}

\address{Paul Kirk: Department of Mathematics, Indiana University, Bloomington, IN 47405 }
\email{pkirk@indiana.edu}

\address{Charles Livingston: Department of Mathematics, Indiana University, Bloomington, IN 47405 }
\email{livingst@indiana.edu}

%%%%%%%ABSTRACT%%%%%%%%%%%%%%

  \begin{abstract}  
We present the complete classification of the subgroup of the classical knot concordance group generated by knots with eight or fewer  crossings.  Proofs are presented in summary.  We also describe extensions of this work to the case of nine crossing knots.
    \end{abstract}

\maketitle
 %%%%%%%SECTION%%%%%%%%%%%%%%
 %%%%%%%SECTION%%%%%%%%%%%%%%

\section{Introduction.}\label{sectionintroduction}

Recall that knots $K$ and $ J $ in $ S^3$ are called {\it concordant} if there is a  smooth embedded disk in $B^4$ with boundary $K\#-J$.  The set of equivalence classes under this relation forms the smooth concordance group $\calc$,   an abelian group with addition  induced by connected sums.  Initial work by Fox and Milnor~\cite{fox-milnor} and Musasugi~\cite{murasugi} developed obstructions that were sufficient to prove that the figure eight knot, $4_1$, represents a nontrivial element of order two in $\calc$ and that the trefoil, $3_1$, is of infinite order.  

  The analysis of the information carried by the universal abelian cover of a knot complement (called the {\em algebraic} information) culminated with Levine's definition ~\cite{levine}    of a homomorphism $\phi$ from $\calc$ onto the {\em  algebraic concordance group} $\calg$, a group isomorphic to $  \zz_2^\infty \oplus  \zz_4 ^\infty \oplus \zz^\infty$.   The nontriviality of the kernel of $\phi$ was proved by Casson and Gordon~\cite{casson-gordon}.

Since then, techniques of increasing   effectiveness have been developed to study $\calc$. For the most part,  these have been applied to consider individual knots or to study specified families of knots, for instance two-bridge knots~\cite{casson-gordon}, pretzel knots~\cite{fs, gj}, positive knots~\cite{rudolph}, Whitehead doubles~\cite{cg, cl, hk}.  Beyond these examples, throughout the study of concordance investigators have built families of knots specifically designed to   realize newly discovered invariants; two examples among many are~\cite{cot, gilm}.      Despite these  many advances, the classification problem for $\calc$ remains far from reach.

Our goal is to  approach  the general classification by completing the classification for the  subgroup generated by eight crossing knots.  In doing so, we will illustrate the effectiveness of twisted Alexander polynomials to resolve previously unmanageable examples of potential concordances.   In the appendix we will summarize   results of the first author~\cite{collins} for the much larger group generated by prime knots of nine and fewer crossing knots.

One might ask precisely what is meant by ``classification.''  From the perspective of fully understanding knot concordance,  at the very least this needs to include   an algorithm to determine if a given linear combination of low crossing number knots is trivial.  But in addition,   the classification should also identify the order of an element and its value in the algebraic concordance group $\calg$ under Levine's homomorphism $\phi$.  Here is a summary statement.  The full result, filling in details concerning orientation issues, will be provided in the next section. 

\begin{theorem}
There are 36 oriented prime knots of eight or fewer crossings.  These knots generate a subgroup $\calc_8 \subset \calc$ isomorphic to $ (\zz^{18} \oplus \zz_2^{6}) \oplus  (\zz^5 \oplus \zz_2)$.   The first of the two summands, $ (\zz^{18} \oplus \zz_2^{6})$,  maps injectively into $\calg$; for the second summand, the $\zz^5$ has image in  $\calg$ isomorphic to $\zz_4 \oplus \zz_2^3$ and the $\zz_2$ maps trivially to $\calg$.
\end{theorem}

The proof of this theorem is such that determining the concordance class of any linear combination of knots, each of eight or fewer crossings, is completely algorithmic.  In fact, an online calculator  developed by the first author, available at {\it KnotInfo}~\cite{cha-livingston}, permits one to easily determine the concordance properties, including the  image in the algebraic concordance group $\calg$,   of any given linear combination of prime knots with eight or fewer crossings.

\section{Statement of classification}

There are 36 prime knots eight or fewer crossings; every knot of eight or fewer crossings can be expressed as a connected sum of these prime knots.  In this count, we include only one knot from each pair, $\{K, \text{mirror}(K)\}$  for those knots $K$ which are reversible,  since the mirror image  with its orientation reversed represents $-K$ in $\calc$.  However, there is one nonreversible knot, $8_{17}$, so for completeness  we must include $8_{17}$ and $8_{17}^r$ in the count.  
 Thus, the   set of knots we consider is  $$3_1, 4_1, 5_{\{1-2\}}, \\ 6_{\{1-3\}},  7_{\{1-7\}}, 8_{\{1-21\}}, 8_{17}^r.$$   

 Let $\calf$ denote the free abelian group generated by these 36 knots.    Levine provided a classification result for $\calg$ based on a set of explicit invariants (enhanced slightly in~\cite{livingston,morita} to consider four torsion).   Using this, an initial decomposition becomes available.
 
 \begin{lemma} There is a direct sum decomposition $\calf  = \calf_\infty \oplus \calf_4 \oplus  \calf_ 2 \oplus  \calf_1$, with summands of rank  $18, 1, 9, $ and $8$, respectively.  The group $\calf_\infty$ maps to a summand of $\calg$ isomorphic to $\zz^{18}$; $\calf_4$ and $\calf_2$ map onto complementary  summands of  $\calg$ isomorphic to $\zz_4$ and $\zz_2^9$, respectively; $\calf_1$ maps trivially to $\calg$.
 
 \end{lemma}

 The next three  results will be seen to be   consequences of Casson-Gordon theory.
 
  \begin{lemma} The group $\calf_4\cong \zz$ injects into $ \calc$.
 \end{lemma}
 
 \begin{lemma}The group $\calf_2$ has a further decomposition as $\calf_2^\infty \oplus \calf_2^2$, of ranks $3$ and $6$, respectively.  The group $\calf_2^\infty$  maps  onto $\zz^3 \subset \calc$  and $\calf_2^2$ maps  surjectively  to $\zz_2^6 \subset \calc$.
 \end{lemma}

 \begin{lemma} The group $\calf_1$ has a further decomposition as $\calf_1^\infty \oplus \calf_1^2 \oplus \calf_1^1$, of ranks $1$, $1$ and $6$, respectively.  The group 
 $\calf_1^\infty$ maps isomorphically onto $\zz  \subset \calc$,   $\calf_1^2$ maps onto   $\zz_2\subset \calc$, and $\calf_1^1$ maps trivially to $\calc$.
 
 \end{lemma}
 
 With these groups defined, the classification is essentially completed by describing the bases of each of these summands.

 \begin{theorem}  \label{main result}  Bases for the subgroups described in the previous lemmas are given in the following list.  (For clarity and consistency we write $\calf_\infty^\infty $ for $\calf_\infty$ and $\calf_4^\infty$ for $\calf_4$.  In this way, subscripts indicate orders in $\calg$ and superscripts give the order in $\calc$.)
 \begin{itemize}
 
 \item {  $\calf_\infty^\infty$\rm{:} }   $3_1, 5_{\{1,2\}}, 6_2, 7_{\{1-6\}}, 8_{\{2,4,5,6,7,14,16,19\}}$\vskip.05in
 
  \item { $\calf_4^\infty$ }\rm{:}   $ 7_7$   \vskip.05in

 \item {$\calf_2^\infty$\rm{:}} $8_1, 8_{13}, (8_{15} - 7_2 - 3_1)$\vskip.05in

 \item {$\calf_2^2$\rm{:}}   $ 4_1, 6_3, 8_{\{3, 12, 17, 18\}}$\vskip.05in
 
 \item {$\calf_1^\infty$\rm{:}}  $(8_{21} - 8_{18} - 3_1)$\vskip.05in

 \item {$\calf_1^2$\rm{:}}  $(8_{17} - 8_{17}^r$)\vskip.05in
 
 \item {$\calf_1^1$\rm{:}} $6_1, 8_{\{8,9, 20\}}, (8_{10} +3_1) , (8_{11} - 3_1)$\vskip.1in
 
 \end{itemize}
 
 \end{theorem}

 \section{Examples.}
Two examples  illustrate how the previous theorem provides a complete classification.   

\subsection{$K = 8_{10} \# 8_{21}$}  Consider the knot $K = 8_{10} \# 8_{21}$.  Writing $K$ in terms of the bases given above, we have $$K = (8_{10} \# 3_1) \# (8_{21} \# -8_{18} \# -3_1) \# 8_{18},$$ the sum of generators of $\calf_1^1$, $\calf_1^\infty$ and $\calf_2^2$.  Thus, $K$ maps to an element of torsion
 two in $\calg$ and is of infinite order in $\calc$.   \vskip.05in

\noindent{\it Remark}:  Both  $8_{10}$ and $8_{21}$   have nontrivial signature functions, identical to those of $3_1$ and $-3_1$, and so individually have infinite order in $\calc$.  However, the sum has trivial signature function, and thus represents a torsion element in $\calg$.  Their Alexander polynomials are the same as those of $3_1 \# 3_1 \# 3_1$ and $3_1 \# 4_1$, respectively.  Thus,   $8_{10}$ and $   8_{21}$ are distinct in $\calc$, since the product of their Alexander polynomials      is not of the form $f(t)f(t^{-1})$ for some $f$ (the Fox-Milnor~\cite{fox-milnor}  obstruction to being slice).     Proving that this knot is of infinite order in $\calc$ is more challenging; the only proof we know depends on a careful analysis of branched covers and twisted Alexander polynomials.  
  \vskip.05in

 \subsection{$K = 8_{17}^r \# 8_{21} \# -3_1$}  In this case, we rewrite $K$ as $$K= - (8_{17} - 8_{17}^r)\  \# \  8_{17}\  \#\ (8_{21} \#-8_{18} \# -3_1)\ \#\  8_{18}.$$  These four knots are in $\calf_1^2$, $\calf_2^2$, $\calf_1^\infty$, and $\calf_2^2$ respectively.  From this it follows that $K$ is of algebraic order two, but is of infinite order in concordance.
 
 \subsection{Calculator}  At these examples illustrate, the problem of determine the (algebraic) concordance class   of a knot and its order has been reduced to simple linear algebra. A program, accessible via the internet, has been written to implement the algorithm.    The reader is invited to visit the ``Concordance Calculator'' posted on {\it KnotInfo} ~\cite{cha-livingston} to generate more examples.
 
 \section{Algebraic classification}  The definitions of $\calf_\infty$, $\calf_4$, $\calf_2$, and $\calf_1$ are algebraic, depending only on the image of the sets of knots in the algebraic concordance group $\calg$.  Levine's classification of $\calg$ is for the most part algorithmic and thus the problem  of identifying the image of $\calf$ in $\calg$ follows from  his    work.
 
 We recall briefly one of the key features of Levine's classification.    The group $\calg$ is defined as certain equivalence classes of Seifert matrices $V_K$ with block sum as the addition.  
 For any field $\ff$ containing $\qq$, one defines a group $\calg_\ff$ which $\calg$ maps into.  For $\ff = \qq$, $\calg_\ff = \calg$.  There is a decomposition $\calg_\ff \cong \oplus_{f(t)} \calg_{f(t)}$ where the $f(t)$ are taken from the set of distinct  $\ff$--irreducible symmetric Alexander polynomials.  A knot maps to summands corresponding to factors of its Alexander polynomial.

For $\ff$ the real numbers, $\rr$, the relevant polynomials $f$  are the irreducible quadratics having roots on the unit circle.  Each of the corresponding summands is isomorphic to $\zz$.     To  a knot $K$  with Seifert matrix $V_K$ and a polynomial $f$ having its roots at $e^{\pm 2\pi i \theta}$, one assigns the integer  which is the jump in the signature function   $$\sigma_\omega(K) = \sigma_{\omega}(V_K)=\text{signature} \left ( (1- \omega )V_K +(1-\omega^{-1})V_K^t\right)$$ as   $\omega $ moves through $ e^{2\pi i \theta} $ on the unit circle. 
The intersection of the kernels of the collection of these integer valued invariants over all such quadratics $f$   
  is precisely the torsion in $\calg$.    This kernel can also be described as the subgroup of $\calg$ generated by those $V_K$ so that $\sigma_\omega(V_K)$ is  zero for all but finitely many $\omega$.  The full set of invariants that detect the two and four torsion arise from considering $p$--adic completions of $\qq$ as the field $\ff$.  For low crossing number knots, this $p$--adic analysis can largely  be avoided.

 \subsection{Slice knots:  $\calf_1^1$.}
 The first step in  producing the decomposition of Theorem \ref{main result}  is to isolate the known {\em slice knots}, that is,   knots that represent 0 in $\calc$.  The first four elements of $\calf_1^1$ can be quickly seen to be slice.  The triviality of $8_{10} \# 3_1$ and $8_{11}\# -3_1$ in $\calc$ was observed by Conway in~\cite{conway}.  It is a consequence of the results described here that there are no other such relations.

 \subsection{Infinite order elements and the signature function: $\calf_\infty^\infty$.}
 
 Since jumps of the signature function occur only at roots of the Alexander polynomial, for any finite set of knots, one can determine the image of the signature function by evaluating it at a set of points that includes a number on the unit circle  between any two such roots of any of the Alexander polynomials that arise.  For the current calculation, and that for nine crossing knots, there are only 70 roots of Alexander polynomials on the upper unit half circle.   Hence, evaluating the signature function at the midpoints of the circular segments joining those roots gives a homomorphism  from $\calf$ to $\zz^{70}$ whose kernel maps by $\phi$  to two and four torsion in $\calg$.

It is now an exercise in elementary linear algebra to identify the set $\calf_\infty^\infty$ as mapping to a generating set of the free part of the image of $\calc$ in $\calg$ and that the remaining summands map to torsion in $\calg$.

\subsection{Torsion in $\calg$}

We now describe the image of the subgroup generated by  $\calf_4^\infty, \calf_2^\infty, \calf_2^2, \calf_1^\infty$ and $\calf_1^2$ in the torsion of  $\calg$.  Levine's paper~\cite{levine} provided a complete set of invariants of algebraic concordance.   Those invariants ranged over all $p$--adic completions of $\qq$, but in~\cite{morita} it was shown how to reduce the set of primes to a finite collection.  In~\cite{livingston} there was a further reduction, leading to effective means of determining the algebraic order of a knot.  The knots we are considering  fall to those techniques, so we only summarize the required work here.  

\subsubsection{$\calf_1^2$, $\calf_1^\infty$}  We most quickly dispense with $\calf_1^2$ and $\calf_1^\infty$.  Since a knot and its reverse represent the same element in $\calg$, $8_{17} - 8_{17}^r$ is  {\em algebraically slice}, that is, maps to zero in $\calg$.     (One proof of this follows from a theorem of Long~\cite{long}.  For any knot $K$, the knot $K \# -K^r$ is ``positive amphicheiral,'' meaning there is an orientation reversing involution of $S^3$ inducing an orientation preserving homeomorphism of $K \#-K^r$.  According to~\cite{long}, such knots are algebraically slice.)   The analysis for $8_{21} - 8_{18} - 3_1$ is a bit more detailed.  However, the Alexander polynomials of   $8_{21}, 8_{18}$, and $ 3_1$    factor into irreducible quadratics, and for these Levine's criteria simplify, as described in~\cite{levine},   from which one can readily show that   $8_{21} - 8_{18} - 3_1$ is algebraically slice.

\subsubsection{$\calf_4^\infty$, $\calf_2^\infty$, $\calf_2^2$} The group $\calg$ splits as a direct sum $\oplus_{p(t)} \calg_{p(t)}$ where the $p(t)$ are distinct irreducible symmetric rational polynomials;  a knot maps to summands corresponding to factors of its Alexander polynomial.   Here is a list of the Alexander polynomials that arise.
\begin{itemize}
\item  $7_7$:   $  1-5t+9t^2-5t^3+t^4$\vskip.05in
\item $8_1$: $  3-7t+3t^2$\vskip.05in
\item $8_{13}$: $    2-7t+11t^2 -7t^2 +2t^4$ \vskip.05in
\item $8_{15} -7_2 - 3_1$:  $   (1-t+t^2)^2(3-5t+3t^2)^2$\vskip.05in
\item  $4_1$: $  1 - 3t +t^2$\vskip.05in
\item $6_3$: $ 1-3t+5t^2-3t^3+t^4$\vskip.05in
\item $8_3$: $ 4-9t +4t^2$\vskip.05in
\item $8_{12}$: $  1-7t+13t^2-7t^3 +t^4$\vskip.05in
\item $8_{17}$: $ 1-4t+8t^2 -11t^3 +8t^4 - 4t^5+t^6$\vskip.05in
\item $8_{18}$: $ 1-5t+10t^2 -13t^3 +10t^4 - 5t^5+t^6$\vskip.05in

\end{itemize}

Since each of these knots maps to a distinct summand of $\calg$, we need only consider the order of the image of each individually.  A few quick observations simplify the work.   First, with the exception of $8_{15} - 7_2 -3_1$, all the polynomials are irreducible, and in particular none are of the form $f(t)f(t^{-1})$, so each maps nontrivially to $\calg$.   Second, the knots in $\calf_2^2$ are all (negative) amphicheiral and so are of order at most  two in $\calc$,    and since they map nontrivially  to $\calg$, they have order two   in $\calc$ as well. 
 
It remains to consider the first four knots on the list.  According to Levine (see~\cite{levine}  for details), four torsion is detected by primes congruent to $3$ mod 4 dividing $|\Delta(-1)|$.  The precise conditions immediately show that $7_7$ is of order four, while $8_1$ and $8_{13}$ are of order two in $\calg$.   The last case, the knot $8_{15} - 7_2 - 3_1$, is the most subtle.  Knots with these Alexander polynomials can represent four torsion in $\calg$, but since the polynomial exponents for this knot  are even, the knot must map either trivially, or to two torsion.  The nontriviality is detected working with the prime 11.

\section{Algebraically slice knots and twisted Alexander polynomials}

It remains to understand the algebraically slice knots represented by elements in $\calf$.  The subgroup of such knots  is generated by $4\calf_4, 2\calf_2,$ and $\calf_1$.  We have already observed that $\calf_1^1$ and $2\calf_2^2$ represent trivial elements in $\calc$.  Thus, we are reduced to considering $4\calf_4^\infty$, $2\calf_2^\infty$, $\calf_1^\infty$ and $\calf_1^2$.  Here is the list of knots, with names now attached for later discussions.

\begin{itemize}
\item  $4\calf_4^\infty$:  $K_1= 4(7_7)$\vskip.05in
\item $2\calf_2^\infty$: $K_2 = 2(8_1), K_3= 2(8_{13}), K_4 = 2(8_{15} - 7_2 - 3_1)$\vskip.05in

\item $\calf_1^\infty $ : $K_5=  8_{21} - 8_{18} - 3_1  $ \vskip.05in

\item $\calf_1^2$: $K_6 =8_{17} - 8_{17}^r$\vskip.05in

\end{itemize}
  
 Individual knots in this subgroup have been discussed in previous articles:   that $7_7$ is of infinite order in $\calc$ is a consequence of Casson-Gordon theory, as presented in~\cite{livingston-naik}; the nontriviality of $8_{17}-8_{17}^r$ was among the first applications of the twisted Alexander polynomial viewed  as a discriminant of a Casson-Gordon invariant, in~\cite{kirk-livingston2}.  Twisted polynomials were used in~\cite{tamulis} to prove that $2(8_{13})$ is not slice, and further Casson-Gordon techniques were applied to show that $8_1$ is of infinite order.   
 
 \subsection{Twisted Alexander polynomials as slicing obstructions}
 
 Let $M_q(K)$ denote the $q$--fold branched cover of $S^3$ branched over $K$, where $q$ is a prime power.  For each prime $p$ and homomorphism $\rho : H_1(M_q(K)) \to \zz_p$, there is a twisted Alexander polynomial,  $\Delta_{K,\rho}(t) \in \qq[\zeta_p][t, t^{-1}]$.  In~\cite{kirk-livingston1} this polynomial is shown to be related to a Casson-Gordon invariant of $K$.  In particular, if $K$ is slice then for certain $\rho$, it is shown that $\Delta_{K,\rho}(t)$ factors as $f(t)\overline{f(t^{-1})}$ for some $f(t) \in  a \qq[\zeta_p][t, t^{-1}]$,where $a$ is some unit in  $\qq[\zeta_p][t, t^{-1}]$.  
 
 Restricting the set of characters that one must consider is one of the most challenging aspects of the computations.  In brief, there is a linking form on the first homology group of the branched cover and there is also a $\zz_q$  action on this group.  The appropriate set of characters is given by linking with elements  in a subgroup of the  homology of the cover, and this subgroup should be invariant under the group action and be self-annihilating with respect to the linking form.  Details can be found in~\cite{kirk-livingston2}.
 
 The computation of $\Delta_{K,\rho}$ is completely algorithmic; in~\cite{herald-kirk-livingston} methods were developed that ensured the rapid   computer calculation    of these polynomials, sufficient to work with all prime knots of 12 or fewer crossings.  For any individual knot, one can often enumerate the possible $\rho$ that must be considered and carry out the computation.  The problem of showing linear independence was considered in~\cite{kirk-livingston3}, where the methods were developed that could isolate properties of the necessary families of characters $\rho$; these methods were applied to families of knots specifically designed to be accessible the new techniques.
 
 \subsection{The homology of the covers}
 
 Fortunately, for the knots of interest, a large number of primes occur in considering only the 2--fold and 3--fold branched covers.  Here are the needed results.
 
 \begin{itemize}
 \item  $|H_1( M_2(7_7)) |= 3 \cdot 7$ \hskip.1in %$|H_1( M_2(7_7))| = 13^2$ 
 
  \item  $|H_1( M_2(8_1))| = 13 $ \hskip.1in %$|H_1( M_2(7_7)) |=2^2 \cdot 5^2$ 

 \item  $|H_1( M_2(8_{13}))| = 29 $ \hskip.1in %$|H_1( M_2(8_{13})) |=2^8 $ 

 \item  $|H_1( M_2(8_{15}))| = 3\cdot 11  $ \hskip.1in %$|H_1( M_2(8_{15} )) |=2^8 \cdot 5^2$ 

\item  $|H_1( M_2(7_{2}))| =  11 $ \hskip.1in %$|H_1( M_2(7_{2} )) |=2^6  $

\item  $|H_1( M_2(3_{1}))| =  3 $ \hskip.1in %$|H_1( M_2(3_{1} )) |=2^2  $

 \item  $|H_1( M_2(8_{21}))| = 3 \cdot 5   $ %\hskip.1in  $|H_1( M_2(8_{21} )) |=2^6 $ 

 \item  $|H_1( M_2(8_{18 }))| = 3^2 \cdot 5   $ %\hskip.1in  $|H_1( M_2(8_{21} )) |=2^8 $ 

 \item  $|H_1( M_2(8_{17 }))| = 37  $ \hskip.1in   $|H_1( M_3(8_{17} )) |=13^2 $ 

 \end{itemize}
 
 \subsection{Summary calculation}
 
Suppose a linear combination of the knots $a_1K_1 + a_2K_2 + a_3K_3 +a_4K_4 + a_5K_5 +a_6K_6  = 0 \in \calc$.  Working with the 2--fold covers and the primes $7, 13,$ and $29$ lets one conclude that the coefficients $a_1$, $a_2$ and $a_3$ all are zero.  Because $K_6 = 8_{17} - 8_{17}^r$, it turns out that working with the prime 37 does not imply that $a_6 = 0$.  However, this is the only knot in the list for which 13 appears in the 3--fold cover, and as described in~\cite{kirk-livingston2}, this is sufficient to show that the coefficient $a_6 = 0$.  (Actually, since $K_6$ is of order two, one shows only that $a_6 $ is even.)

At this point, one need only consider combinations $a_4 K_4 + a_5 K_5$.  The primes $11$ and $5$ let one reduce this to considering each knot individually, for which another computation of twisted Alexander polynomials, along with detailed work at identifying possible metabolizers, completes the project.

  \appendix
  \section{Nine crossing knots}
  
  The analysis of the subgroup of $\calc$ generated by prime knots of nine  or fewer crossings   was undertaken in the  Ph.D.  thesis of the first author, available at~\cite{collins}.  In this case there are 87 knots to consider, including three arising from nonreversible knots: $8_{17}^r, 9_{32}^r$ and $9_{33}^r$.    There   are again seven relevant subgroups  of the free abelian group $\calf$ generated by these 87 knots (with superscripts indicating their order in $\calc$ and subscripts the order in $\calg$):

\begin{itemize}

\item $\calf_\infty^\infty$: $3_1, 5_{1-2}, 6_2, 7_{1-6}, 8_2, 8_{4-7}, 8_{14}, 8_{16}, 8_{19}, 9_{1-2} , 9_{4-7}, 9_{9-11},
  9_{13}, 9_{15},\\  9_{17-18}, 9_{20-22}, 9_{25-26}, 9_{31-32}, 9_{35-36}, 9_{38}, 9_{43}, 9_{45}, 9_{47-49}$  \vskip.05in

\item $\calf_4^\infty$: $7_7, 9_{34}$  \vskip.05in

\item $\calf_2^\infty$: $8_1, 8_{13}, (8_{15} - 7_2 - 3_1), {\bf (9_2 - 7_4)}, (9_{12} - 5_2), 9_{14}, (9_{16} - 7_3 - 3_1), \\ 9_{19}, (9_{28} - 3_1), 9_{30}, 9_{33}, (9_{42} + 8_5 - 3_1), (9_{44} - 4_1)$  \vskip.05in

\item $\calf_2^2$: $4_1, 6_3, 8_3, 8_{12}, 8_{17-18} $  \vskip.05in

\item $\calf_1^\infty $: ${\bf (8_{21} - 8_{18} - 3_1)}, (9_8 - 8_{14}), {\bf (9_{23} - 9_2 - 3_1)}, (9_{29} - 9_{28} +2(3_1)),\\{\bf (9_{32}^r - 9_{32})},  (9_{33}^r - 9_{33}), (9_{39} + 7_2 - 4_1), {\bf (9_{40} - 8_{18} - 4_1 - 3_1) }$  \vskip.05in

\item $\calf_1^2 $: $(8_{17}^r - 8_{17})$  \vskip.05in

\item $\calf_1^1 $: $	6_1, 8_{8-9}, (8_{10} + 3_1), (8_{11} - 3_1), 8_{20}, (9_{24} - 4_1), 9_{27}, (9_{37} - 4_1), 9_{41}, 9_{46}$  \vskip.05in

\end{itemize}

  The classification is  nearly  complete,  
  with the exception that the span in $\calc$ of the five knots shown in bold is not fully identified.   Of these, the first is of algebraic order two and the rest are algebraically slice.  Each of the individual knots is of infinite order in $\calc$ and it is presumed that they are linearly independent.  However, in~\cite{collins} and in all subsequent work there has been no success in ruling out any possible relations between these five, except for the fact that any such relation must contain an even multiple of the first.  
  
   Thus the concordance classification of nine crossing knots will be completed by answering the following question: For what   integers $a_1,a_2,a_3,a_4, a_5$ is the linear combination 
  $$ 2a_1J_1+a_2J_2+a_3J_3+a_4J_4+a_5J_5$$ slice, with 
  $$J_1=9_2 - 7_4, J_2= 8_{21}-8_{18} -3_1, J_3=9_{23}-9_2-3_1,$$
  $$ J_4=9^r_{32}-9_{32}, J_5= 9_{40} - 8_{18} - 4_1 - 3_1?$$

%%%%%%%END%%%%%%%%%%%%%%
%%%%%%%END%%%%%%%%%%%%%%
%%%%%%%END%%%%%%%%%%%%%%
%%%%%%%END%%%%%%%%%%%%%%
%%%%%%%END%%%%%%%%%%%%%%
%%%%%%%END%%%%%%%%%%%%%%

%%%%%%%BIBLIOGRAPHY%%%%%%%%%%%%%%

\newcommand{\etalchar}[1]{$^{#1}$}

\end{document}